\documentclass{article}
\usepackage{amssymb,latexsym,amsfonts,amsmath}
\usepackage[paperwidth=8.5in,paperheight=10.5 in]{geometry}
\usepackage{multicol}
\usepackage[utf8]{inputenc}
\usepackage[T1]{fontenc}
\usepackage{microtype}
\usepackage{color}
\usepackage[dvipsnames]{xcolor}
\usepackage{hyperref}
 \hypersetup{
     colorlinks=true,
     linkcolor=Purple,
     filecolor=black,
     citecolor = black,      
     urlcolor= Purple,
     }
\usepackage{cite}
\usepackage{mdframed} 
\usepackage{ntheorem} 
\theoremstyle{nonumberplain}
\newmdtheoremenv[ntheorem]{thm}{Theorem}
\usepackage{tkz-euclide}
\usepackage{tikz}
\usepackage{array}
\usepackage{pgfplots}
\usetikzlibrary{arrows}

\pgfplotsset{my style/.append style={axis x line=middle, axis y line=
middle, xlabel={$x$}, ylabel={$y$}, axis equal }}

\newmdtheoremenv{fact}[thm]{Fact.}
\newmdtheoremenv{boxz}[thm]{ \ }
\newmdtheoremenv{htp}[thm]{How to Prove}

\newtheorem{prop}[thm]{Proposition}

\newmdtheoremenv{rem}[thm]{Thing You Should Remember.}
\newmdtheoremenv{defn}[thm]{Definition.}
\newmdtheoremenv{inabox}[thm]{ \ }

\newcommand{\be}{\begin{equation}}
\newcommand{\ee}{\end{equation}}
\newcommand{\ba}{\begin{eqnarray}}
\newcommand{\ea}{\end{eqnarray}}
\newcommand{\ban}{\begin{eqnarray*}}
\newcommand{\ean}{\end{eqnarray*}}
\newcommand{\barr}{\begin{array}}
\newcommand{\earr}{\end{array}}

\usepackage{caption}
\newcommand{\N}{{\mathbb N}}

\usepackage{amsmath}
\title{A Mathematical Foundation for the Numberlink Game}

\author{Andrea Arauza Rivera, Matt McClinton, David Smith }

\date{\today}
\usepackage[utf8]{inputenc}
\usepackage[english]{babel}
\usepackage{float}

\begin{document}
\maketitle

\hrule
\medskip

\noindent *Andrea Arauza Rivera, Ph. D. \\
Assistant Professor, Mathematics\\
Cal State East Bay\\
andrea.arauzarivera@csueastbay.edu\\

\noindent Matt McClinton\\
Cal State East Bay\\
mmcclinton2@horizon.csueastbay.edu\\

\noindent David Smith\\
Cal State East Bay\\
david.smith2@csueastbay.edu 
\medskip

\noindent *Corresponding author.\\
\hrule

\begin{abstract}
 Numberlink is a puzzle game in which players are given a grid with nodes marked with a natural number, $n$, and asked to create $n$ connections with neighboring nodes. Connections can only be made with top, bottom, left and right neighbors, and one cannot have more than two connections between any neighboring nodes. 
 In this paper, we give a mathematical formulation of the puzzles via graphs and give some immediate consequences of this formulation.
 The main result of this work is an algorithm which provides insight into characteristics of these puzzles and their solutions.
 Finally, we give a few open questions and further directions. 
\end{abstract}

\tableofcontents

\newpage

\section{Introduction}

With millions of downloads between the Apple and Google Play stores, Puzzledom is an app that offers a variety of curious puzzles. One of these is the puzzle called Numberlink \cite{pzzldom}. Puzzles like Numberlink have been the subject of a number of interesting articles; the reader may enjoy any one of the following \cite{velleman2020bicycle}, \cite{plambeck2020barrycades}. The reader should note that there is another popular game called Numberlink. 
This version of the game is discussed in \cite{adcock2015zig}, \cite{ruangwises2020physical}, \cite{yoshinaka2012finding}. 

\begin{figure}
    \centering
    \includegraphics[scale = .4]{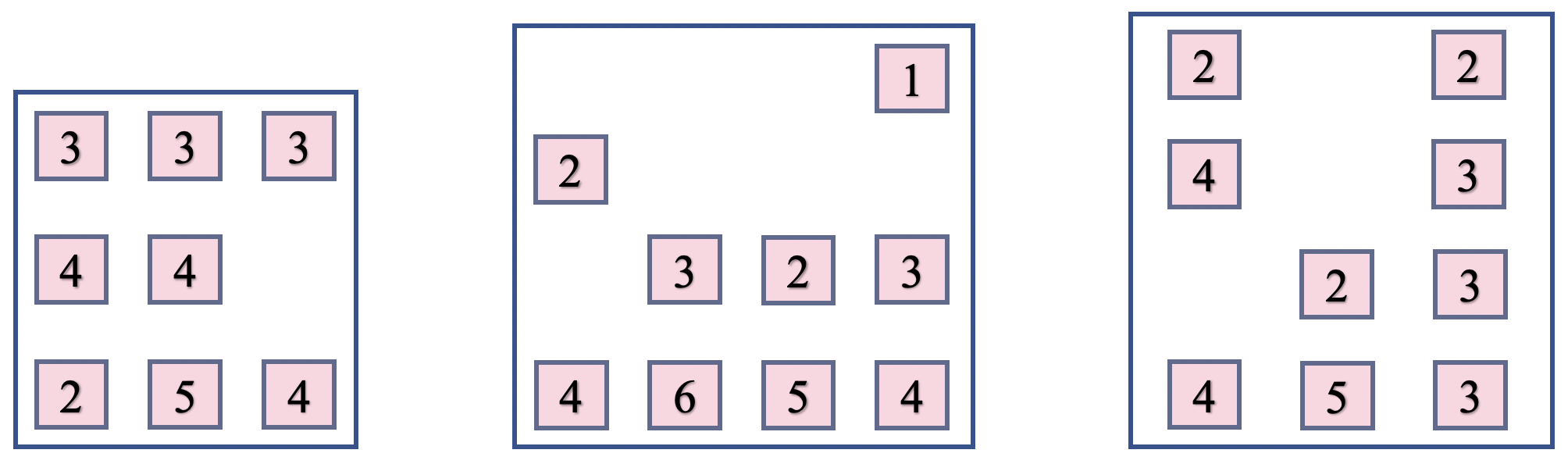}
    \caption{We show 3 sample puzzles for the reader to try and enjoy. Each of these puzzles is included in Numberlink and listed as ``Novice'' (left) or ``Regular''(center and right) \cite{pzzldom}.}
    \label{fig:example2}
\end{figure}

We focus on the Numberlink puzzles found in the Puzzledom app. These puzzles begin with a set of numbered boxes (nodes) configured in a grid; see Figure \ref{fig:example2}. The puzzle is solved when the player creates links between the numbered nodes so that 
\begin{itemize}
    \item a node with number $n$ has $n$ connections to other nodes, 
   \[ \includegraphics[scale = .35]{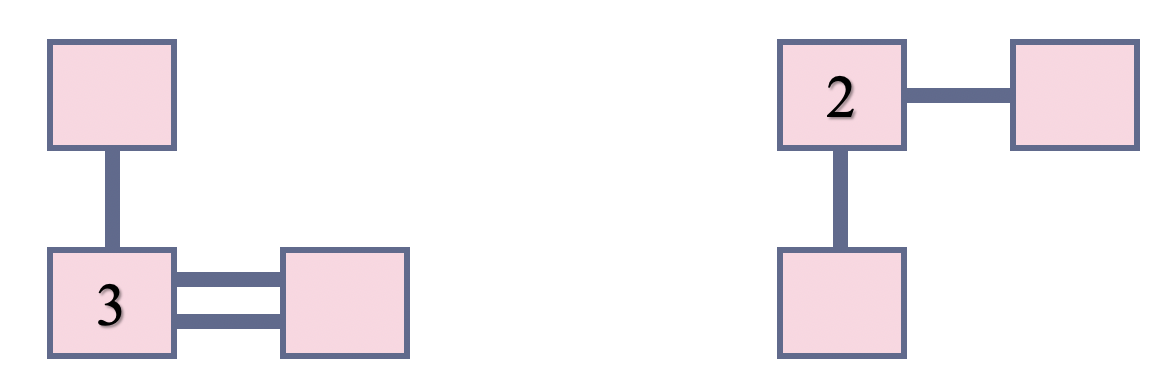}\]
    \item there are no more than $k=2$ connections between any two nodes, and 
   \[ \includegraphics[scale = .35]{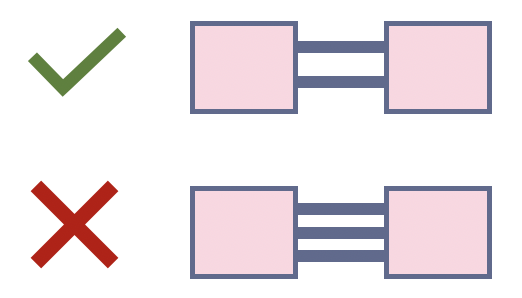}\]
    \item the connections must create a path between any two given nodes (path connected).
   \[ \includegraphics[scale = .35]{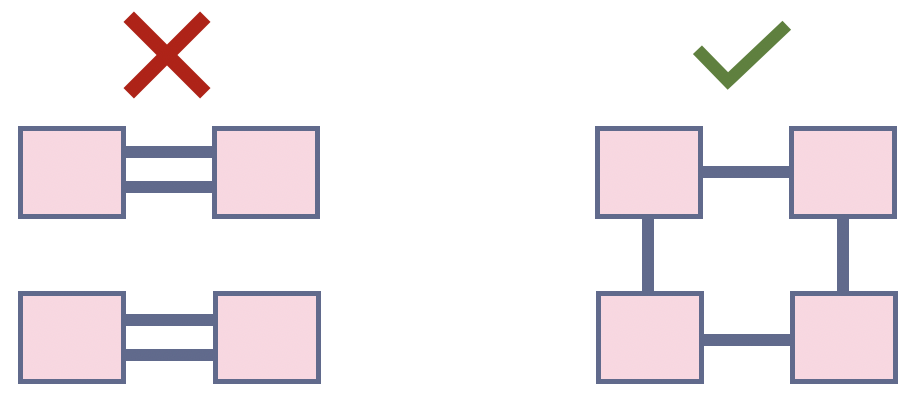}\]
    \item The game also implicitly requires that connections be made horizontally or vertically, and that no connections intersect. 
    \[ \includegraphics[scale = .35]{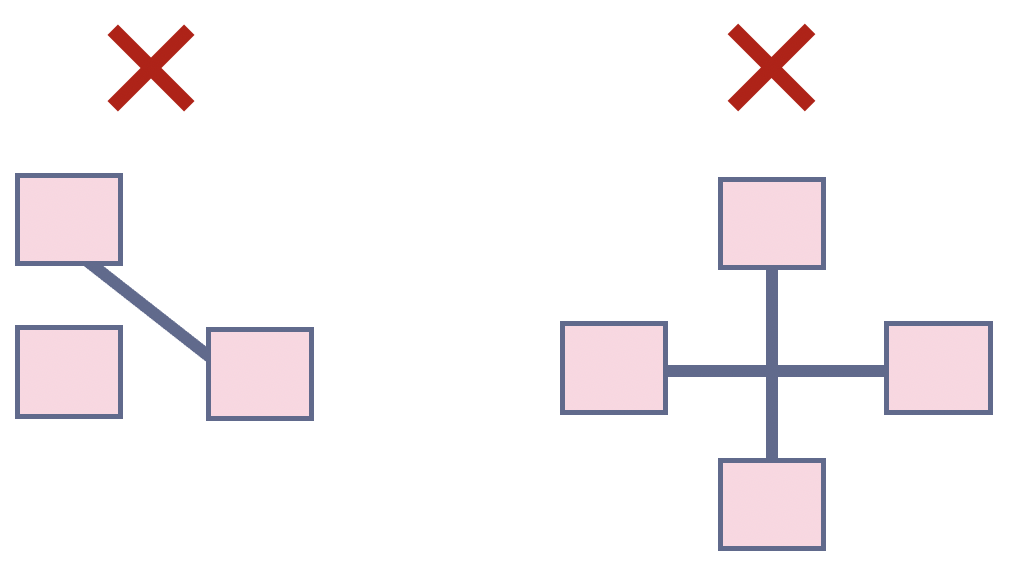}\]
\end{itemize}

Three sample Numberlink puzzles are shown in Figure \ref{fig:example2}. The reader is encouraged to whip out a pencil and try to solve each puzzle! 
The sections in this article are set up as follows:
\begin{itemize}
    \item Section 1 gives a mathematical formulation of the Numberlink game via graphs. This section also includes some initial consequences of this formulation. 
    \item Section 2 describes ways in which the player may find the best nodes to start to create connections. This analysis is based on the number in the node as well as the number of neighbors available for connecting. 
    \item Section 3 contains the main results of this article; an algorithm which produces the ``guaranteed connections'' between nodes. In this section we prove that if the algorithm we outline arrives at a solution, then the solution is unique. 
    \item We conclude in section 4 with some closing thoughts and open questions about these puzzles. 
\end{itemize}

\section{Numbered $k$-Grids}
This section contains the definitions needed to describe the Numberlink puzzles in terms of graphs.  We begin by defining the initial set up of a puzzle as a graph where the nodes are labeled with a whole number $n$ and arranged on a grid. We call these \emph{numbered $k$-grids}. Next, we define what it means for a node $p$ in a numbered $k$-grid to have \emph{top, bottom, left and right neighbors}. Finally, we define what it means for a numbered $k$-grid to be \emph{solved}. 

The reader may be wondering what on earth this $k$ business is. Indeed, in the original Puzzledom-Numberlink puzzles there is a rule that no two nodes may share more than $k=2$ connections. We work with a more general rule and allow $k$ to be any positive whole number.

\begin{defn}
A \textbf{numbered k-grid} is a finite collection of nodes and connections between nodes satisfying the following:
\begin{enumerate}
    \item each node is given coordinates $(x, y)$ where $x, y \in \N \cup \{0\}$;
    \item each node is labeled with a 
    magnitude $n\in \N$;
    \item there are no more than $k$ connections between nodes;
    \item connections can only exist connecting horizontally or vertically adjacent nodes. 
\end{enumerate}
We denote a numbered $k$-grid by $\Gamma_k$ and write $p = [x, y, n]$ for a node in $\Gamma_k$. Often we will denote the magnitude of $p$ as $magn(p)$.
\end{defn}

\begin{figure}
    \centering
    \includegraphics[scale = .4]{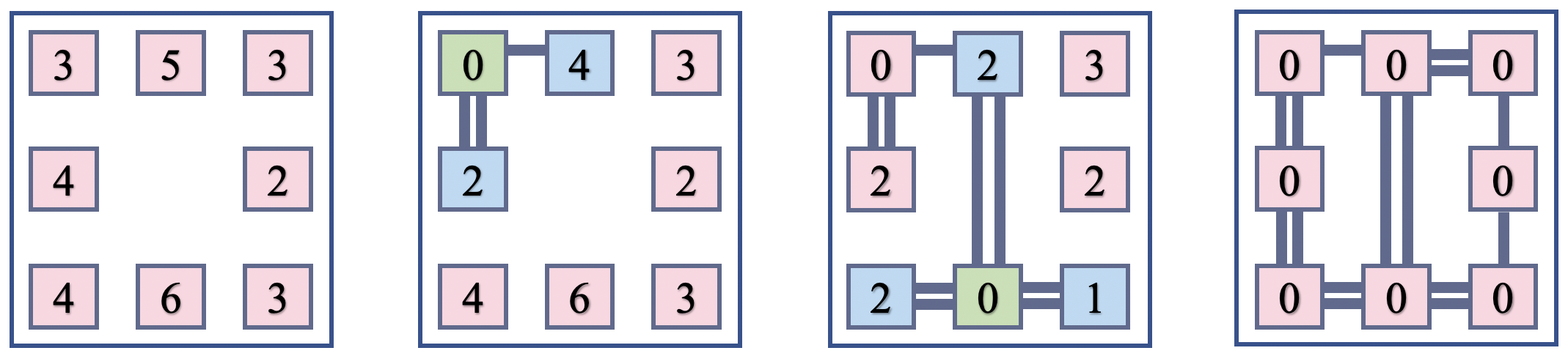}
    \caption{A typical initial set-up for a Numberlink game is shown (left). Two examples of fully connected nodes, given this initial set-up, are shown (center-left and center-right). Finally, a full solution of the game with the initial set-up is shown (right).  }
    \label{fig:example1}
\end{figure}

The choice of a square to represent the nodes in a numbered grid is somewhat arbitrary. One could use a circle or other figure to represent nodes. We choose to use a square to emphasize the top, bottom, left and right neighbors of the node. 
We now define what it means to be neighboring to a node $p$.
\begin{defn}
\textbf{Neighboring Nodes in $\Gamma_k$} \\
Let $W = \{1,2,3,4\}$ and $\Gamma_k$ be a numbered $k$-grid. Consider any node $p = [x_p,y_p,n_p] \in \Gamma_k$. We refer to $1_p,2_p,3_p,4_p$ as the \textbf{top, right, bottom, and left neighbors to p}, respectively. More formally, for each $l \in W$, we define the $l_p$ neighbor by the following:
\begin{align*}
    &1_p = [x_{1_p}, \ y_{1_p}, \ n_{1_p}] \in \Gamma_k \text{ where } x_{1_p}=x_{p}, \text{ and } y_{1_p} = \min \{y : [x_{p}, y, m] \in \Gamma_k \text{ and }y > y_p \}. \\
    & \newline \\
    &2_p = [x_{2_p}, \ y_{2_p}, \ n_{2_p}] \in \Gamma_k \text{ where } y_{2_p}=y_{p}, \text{ and }x_{2_p} = \min \{x : [x, y_p, m] \in \Gamma_k \text{ and }x > x_p \}. \\
    & \newline \\
    &3_p = [x_{3_p}, \ y_{3_p}, \ n_{3_p}] \in \Gamma_k \text{ where } x_{3_p}=x_{p}, \text{ and } y_{3_p} = \max \{y : [x_p, y, m] \in \Gamma_k \text{ and }y < y_p \}. \\
    & \newline \\
    &4_p = [x_{4_p}, \ y_{4_p}, \ n_{4_p}] \in \Gamma_k \text{ where } y_{4_p}=y_{p}, \text{ and } x_{4_p} = \max \{x : [x, y_p, m] \in \Gamma_k \text{ and }x < x_p \}.
\end{align*}
If no node in $\Gamma_k$ satisfies the conditions of an $l_p$ neighbor, we say the $l_p$ neighbor does not exist.
\end{defn}

In Numberlink, when one creates a connection between nodes, the magnitude of each of the two nodes is reduced by 1. A puzzle is completed when the magnitude of all nodes has been reduced to 0 and no rules of the puzzle have been violated (e.g. no nodes share more than $k$ connections, the graph is path connected, connections are only shared by neighboring nodes, and no connections cross). 
We now give the definition of a solved $k$-grid in our context of graphs.
\begin{defn}
\textbf{A solved $k$-grid} \\
For any numbered $k$-grid, $\Gamma_k$, we say it is \textbf{solved} if the set of connections in $\Gamma_k$ satisfy the following:
\begin{itemize}
    \item Every node $p \in \Gamma_k$ is a \textbf{completed node}, meaning the number of connections between $p$ and it's neighboring nodes equals $magn(p)$.
    \item The number of connections between any two neighboring nodes is at most $k$.
    \item No pair of connections 
    intersect.
    \item The graph $\Gamma_k$ is path connected.
\end{itemize}
\end{defn}


\subsection{Some initial results}

\begin{figure}
    \centering
    \includegraphics[scale = .35]{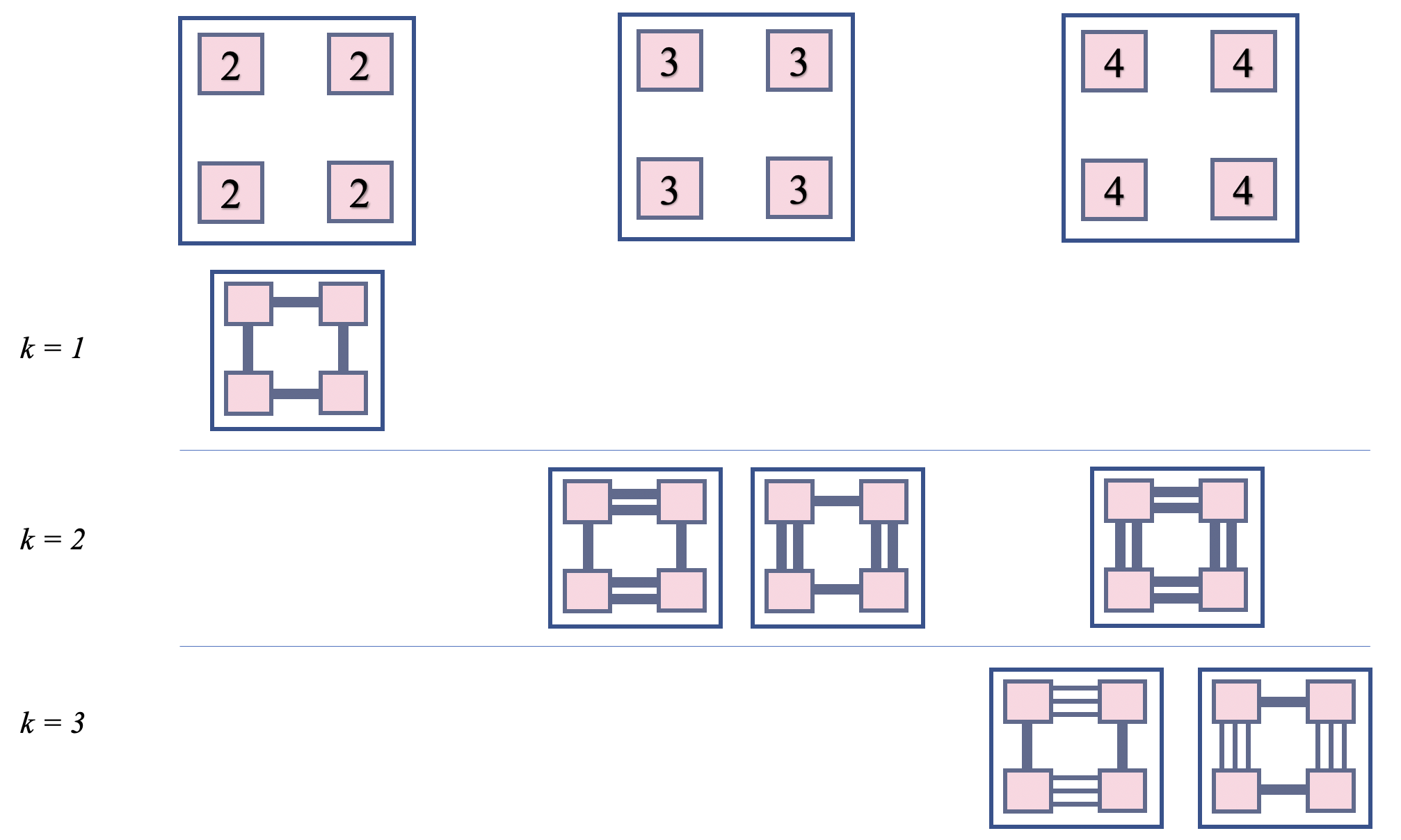}
    \caption{Three examples of the ways in which $k$ affects the number of solutions that can exist for a $k$-grid.}
    \label{fig:exksolutions}
\end{figure}

Every puzzle in Numberlink can be classified as a solvable $2$-grid. 
Note also that any solvable $1$-grid is also a solvable $2$-grid. More generally, if $\Gamma_k$ is solvable, then the graph with the same nodes as $\Gamma_k$ can be solved with any $k'>k$. 

Let's take a moment to inspect how $k$ influences solvability. Consider the 3 examples at the top of Figure \ref{fig:exksolutions}. 
We note the following:
\begin{itemize}
    \item For $k=1$, only the left grid is a solvable $1$-grid. 
    \item Moving to $k=2$ with the left grid, 
    the graph resulting from making 2 connections between nodes would be disjoint. Hence the left grid is a solvable $2$-grid if and only if a single connection is formed between all nodes. 
    \item With $k=2$, the middle and right grids are now solvable. The middle grid has two solutions and the right grid is uniquely solved. 
    \item If $k=3$, we see that the right grid has another 2 solutions. No new solutions are available for the left or middle grids. 
\end{itemize}

The above are examples of $k$-grids whose solutions are relatively simple to build by inspection. As one progresses through the vast array of puzzles in Numberlink, eyeballing solutions involves a good deal of mental planning. Every puzzle in Numberlink has at least one solution, yet an arbitrary $k$-grid has no guarantee of being solvable. Are there ways of quickly identifying unsolvable $k$-grids by inspection? In fact, there are! 

\begin{prop}
Let $\Gamma_k$ be any numbered $k$-grid. If any of the following hold, then $\Gamma_k$ is unsolvable.
\begin{enumerate}
    \item $\Gamma_k$ has a node with no neighbors.
    \item The sum of the magnitudes of all nodes in $\Gamma_k$ is odd.
    \item There is a node $p \in \Gamma_k$ such that the sum of the magnitudes of the neighbors of $p$ is less than the magnitude of $p$.
    \item Every configuration of connections for some node $p$ results in a disjoint graph. 
    \item A node $p$ has $r$ neighbors $(\text{where } r = 1,2,3,4)$ and $magn(p) > rk$.
    \item (Incompatible nodes property) If $k>1$ and there exists a node $p$ with $r$ neighboring nodes such that $magn(p)=(r-1)k+j$ for $j=2,3,\dots,k$, and $p$ has a neighboring node $q$ with $magn(q) \leq j-1$.
\end{enumerate}
\end{prop}

\noindent \textit{\textbf{Proof. }}
For part 1, if a node in any $k$-grid has no neighboring nodes, then it is impossible to form any connections, and thus unsolvable for any $k$. 

Part 2 holds true by the ``first theorem of graph theory'' \cite{merris} 
. Since each connection reduces the magnitudes of the two corresponding nodes by 1, there will be no way to complete every node while also having the sum of the magnitudes be odd. 

For part 3, let $\Gamma_k$ be some $k$-grid with a node $p$ where the sum of the magnitudes of the neighboring nodes is less than the magnitude of $p$, i.e $$magn(p) > \sum\limits_{l \in W} magn(l_p).$$

To say $\Gamma_k$ is solvable means there exist a connection configuration for $p$ where $p$ is complete. Yet the neighboring nodes to $p$ cannot take in $magn(p)$ connections.

The proof for part 4 also follows from the definition of a solved $k$-grid. Let $\Gamma_k$ be a $k$-grid with a node $p$ such that every edge formation results in $\Gamma_k$ disjoint. Meaning all completions of $p$ violates the criteria of being solved and hence $\Gamma_k$ unsolvable. 

Now consider any $k$-grid $\Gamma_k$ satisfying the condition in part 5 for $r=1$. To say $p$ has one neighbor with $magn(p) > k$ implies $magn(p) \geq k+1$. Meaning that there must be more than $k$ connections between $p$ and its neighboring node in $\Gamma_k$, and hence $\Gamma_k$ is unsolvable. For $r=2$, $magn(p) > 2k$ implies $magn(p) \geq 2k+1$. Again this means that more than $k$ connections must be formed for at least one neighboring node, which deems $\Gamma_k$ as unsolvable. Similar arguments can be made for $r=3,4$.


For part 7, if $magn(p) = (r-1)k+j$ with $j = 2,3,4, \dots,k$, then $p$ will need to have $k$ connections with $(r-1)$ of it's neighbors and $j$ connections with the remaining neighbor. Thus if $p$ has a neighbor $q$ with $magn(q) \leq j-1$, then there will be no way to complete the node $p$. 
$\quad \blacksquare$ \\

The reader is encouraged to find their own condition which guarantees that a $k$-grid will be unsolvable. Note that the incompatible nodes property tells us about what magnitudes two neighboring nodes cannot have and still allow the grid to be solvable. For example, in a 2-grid, a node of magnitude 8 with four neighbors cannot have a neighbor of magnitude 1. 
The reader should take a moment to write out other instances in which the incompatible nodes property tells us that a $k$-grid will be unsolvable. 

\subsection{Creating connections between nodes}

Now let's talk about connections. Nodes of relatively high and low magnitudes often have little variety in how the connections are constructed. A node of magnitude 3 with four neighbors in a $2$-grid has only 4 ways of creating connections. However, a node of magnitude 20 with four neighbors in a $10$-grid has 891 ways of drawing connections! It is not always the case that all neighboring nodes can take in a full $k$ connections. This means that for the magnitude 20 node, we don't always have to consider all 891 cases. So how do we avoid exhausting every feasible case? We will now begin to describe an algorithm to determine what connections a node \emph{must} have. 

\begin{defn}
\textbf{Connections between Neighboring Nodes} \\ 
Let $\Gamma_k$ be a $k$-grid and $p \in \Gamma_k$. Let $W = \{1,2,3,4\}$. Define $\Phi_k[p]$ to be the set of all connection configurations for the node $p$ within $\Gamma_k$. 

Thus we have  
$$\Phi_k[p] = \{ \omega = \omega_1 \omega_2 ... \omega_n : \omega_i \in W, n=magn(p), \omega_i = \omega_j \text{ at most k times}\}. $$
Here the symbols $1, 2, 3, 4$  refer to the top, right, bottom or left neighbor of the node $p$, respectively. 
We denote the \textbf{length} of a word by $||\omega||$. We say two words are equivalent, $\omega = \omega'$ in $\Phi_k[p]$, when they have the same symbols without deleting duplicates $\{\omega_i\} = \{\omega'_j\}$. 
\end{defn}

We emphasize that the rearrangement of the $\omega_i$ does not result in a new word. To say $\omega = 1212121$ means that 4 connections are created with the top ($\omega_i = 1$) neighbor, and 3 connections are created with the right ($\omega_i = 2$) neighbor. This reflects the idea that the order in which the connections are drawn is irrelevant. 

When considering the ways in which we create connections between nodes, we look to a subset of $\Phi_k[p]$ which is aware of not only the magnitude of a node and the value of $k$, but also the ``surroundings'' of the node $p$ and the consequences of making the connections indicated by an $\omega$. 

\begin{defn}
\textbf{Feasible Connection Configurations in $\Gamma_k$}. \\
Define $count(\omega,l)$ for $l \in W$ to be the number of occurrences of symbol $l$ in the word $\omega$. Let $\Gamma_k$ be any $k$-grid and let $p \in \Gamma_k$. The set of all \textbf{feasible connection configurations} is denoted by $\Phi_{\Gamma_k}[p] $ and defined as the set of $\omega \in \Phi_k[p]$ satisfying the following: 

\begin{enumerate}
    \item Creating the connections in $\omega$ will not violate the upper bound on connections, $k$.
    \item Creating the connections in $\omega$ will not create intersecting connections between nodes. 
    \item $count(\omega,l) \leq \min\{k,magn(l_p)\}$ for every $l \in W$.
    \item The $k$-grid resulting from making the connections in $\omega$ in $\Gamma_k$ does not result in a set of 
    solved nodes which are disjoint from other nodes. 
    \item No node is incomplete with complete neighbors after making the connections in $\omega$.
\end{enumerate}
Set $\omega^*(p) = \underset{\omega \in \Phi_{\Gamma_k}[p]}{\bigcap} \{ \omega_i \} $ without deleting duplicates. 
\end{defn}

The conditions required for a word $\omega$ to be in $\Phi_{\Gamma_k}[p]$ follow from the intuition one builds from playing the Numberlink puzzles. The set $\Phi_{\Gamma_k}[p]$ is often a refinement of $\Phi_k[p]$ where now the words in $\Phi_{\Gamma_k}[p]$ obey more of the rules for solving a $k$-grid. Looking at what connections are common to each word in $\Phi_{\Gamma_k}[p]$, gives us a list of connections that are \emph{guaranteed} to be made in a solution for the $k$-grid. 


\begin{figure}
    \centering
    \includegraphics[scale = .4]{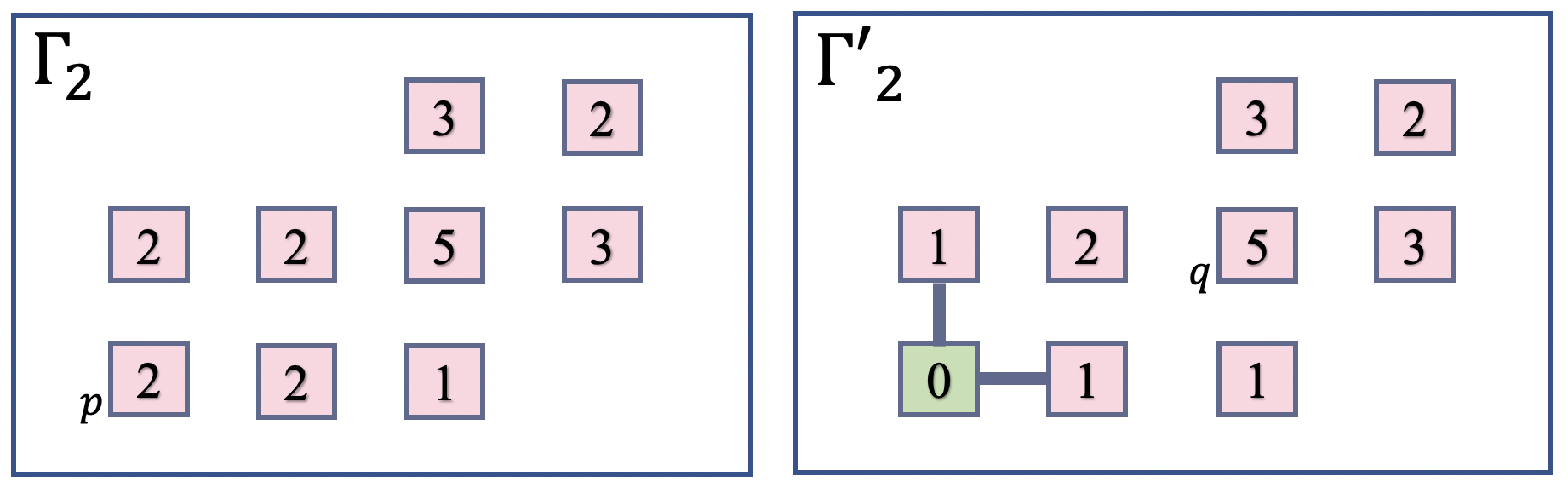}
    \caption{We show a 2-grid, $\Gamma_2$, on the left with node $p$ in the lower bottom. Applying the word 12 to $p$ we get $\Gamma'_2$ on the right. The node $q$ in $\Gamma'_2$ refers to the center node of magnitude 5.}
    \label{fig:defofphi}
\end{figure}
\subsection{Example of computing $\Phi_k[p]$ and $\Phi_{\Gamma_k}[p]$}

Let's do an example that illustrates how one constructs the sets $\Phi_k[p]$ and $\Phi_{\Gamma_k}[p]$. Any node of magnitude 2 in a 2-grid will have
$$\Phi_2[p] = \{11, 22, 33, 44, 12, 13, 14, 23, 24, 34\}.$$ 
Consider the $2$-grid on the left in Figure \ref{fig:defofphi}. 
The node of magnitude 2 in the left bottom corner will have the following as $\Phi_{\Gamma_2}[p]$:
$$\Phi_{\Gamma_2}[p] = \{ 12\}.$$
The words 33, 44, 13, 14, 23, 24, and 34 are not included in the set because of condition 3 in the definition of $\Phi_{\Gamma_2}[p]$ and the words 11 and 22 are not included since they violate condition 4 of the definition of $\Phi_{\Gamma_2}[p]$. We then see that $\omega^*(p)$ for this lower bottom node of magnitude 2 is $\omega^*(p) = 12$. Making these connections gives the $2$-grid on the right of Figure \ref{fig:defofphi}. 

Let's repeat these computations for the center node of magnitude 5, where this time we are working from the 2-grid on the right which includes some connections. We get,
\begin{align*}
    \Phi_2[q] = \{11223,11224,11233,11234,11244,11334, & 11344,12233,12234,  12244, \\ 
    & 12334,12344,13344,22334,22344,23344\}
\end{align*}
and
$$\Phi_{\Gamma'_2}[q] = \{11223, 11224, 11234, 12234\}.$$
In this case, the words that include a 33 will violate condition 3 in the definition of $\Phi_{\Gamma'_2}[q]$ and the words that include a 44 will violate condition 5. 
Here we get that $\omega^*(q) = 12$ for the node $q \in \Gamma_2'$. This tells us that any solution to the grid $\Gamma'_2$ will include the connections 12 for node $q$. The reader is encouraged to find a solution to $\Gamma_2$ in Figure \ref{fig:defofphi}.  




\section{The Path Towards a Solution}


In this section we begin with an examination of where to start forming connections in a $k$-grid. We then move to the main result in this work---an algorithm that finds and creates all guaranteed connections between nodes. 

\subsection{Where do we start?}
Let's review the mechanisms we have so far for computing connections between nodes. Based only on the magnitude of a node and $k$ we can generate the list of words in $\Phi_k[n]$. This gives us any and all possible ways to form $n$ connections with 4 or fewer neighbors while respecting the rule that at most $k$ connections are shared between two nodes. We refine this list by looking at the subset $\Phi_{\Gamma_k}[p]$ which is now aware of the neighbors of $p$ and throws away any connection configurations which would immediately result in some violation of the rules for $k$-grids. Finding what is common to all words in $\Phi_{\Gamma_k}[p]$ gives us $\omega^*(p)$. 

What is missing thus far is an indication of what node to choose at the start! One way to approach this is to find which nodes have the fewest possible configurations of connections given the node's magnitude, the value of $k$, and the number of neighbors that the node has. For example, a node of magnitude 7 with 4 neighbors in a $2$-grid will have only 4 ways to configure the connections. The large magnitude paired with the small value of $k$ means that there is only 1 connection that is different for each configuration. 
\[\includegraphics[scale = .4]{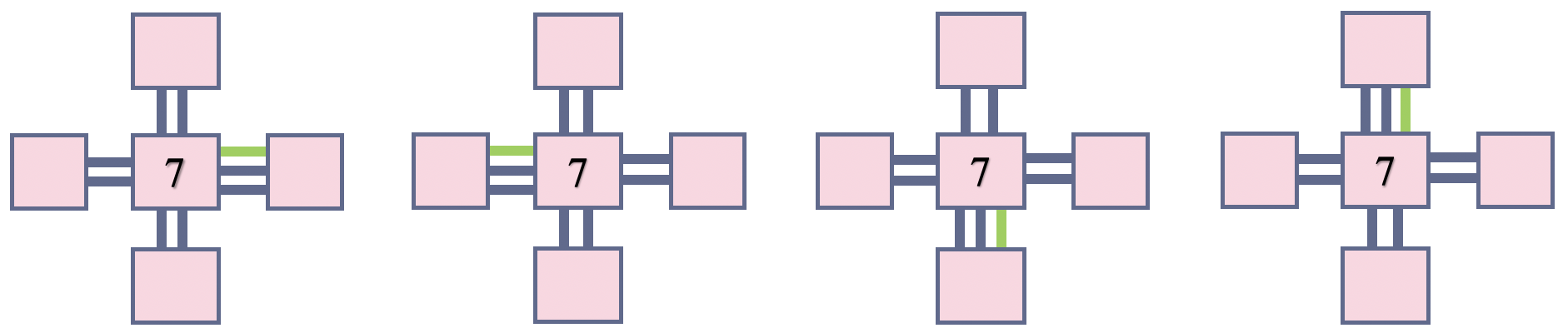}\]

\begin{figure}
    \centering
    \includegraphics[scale = .4]{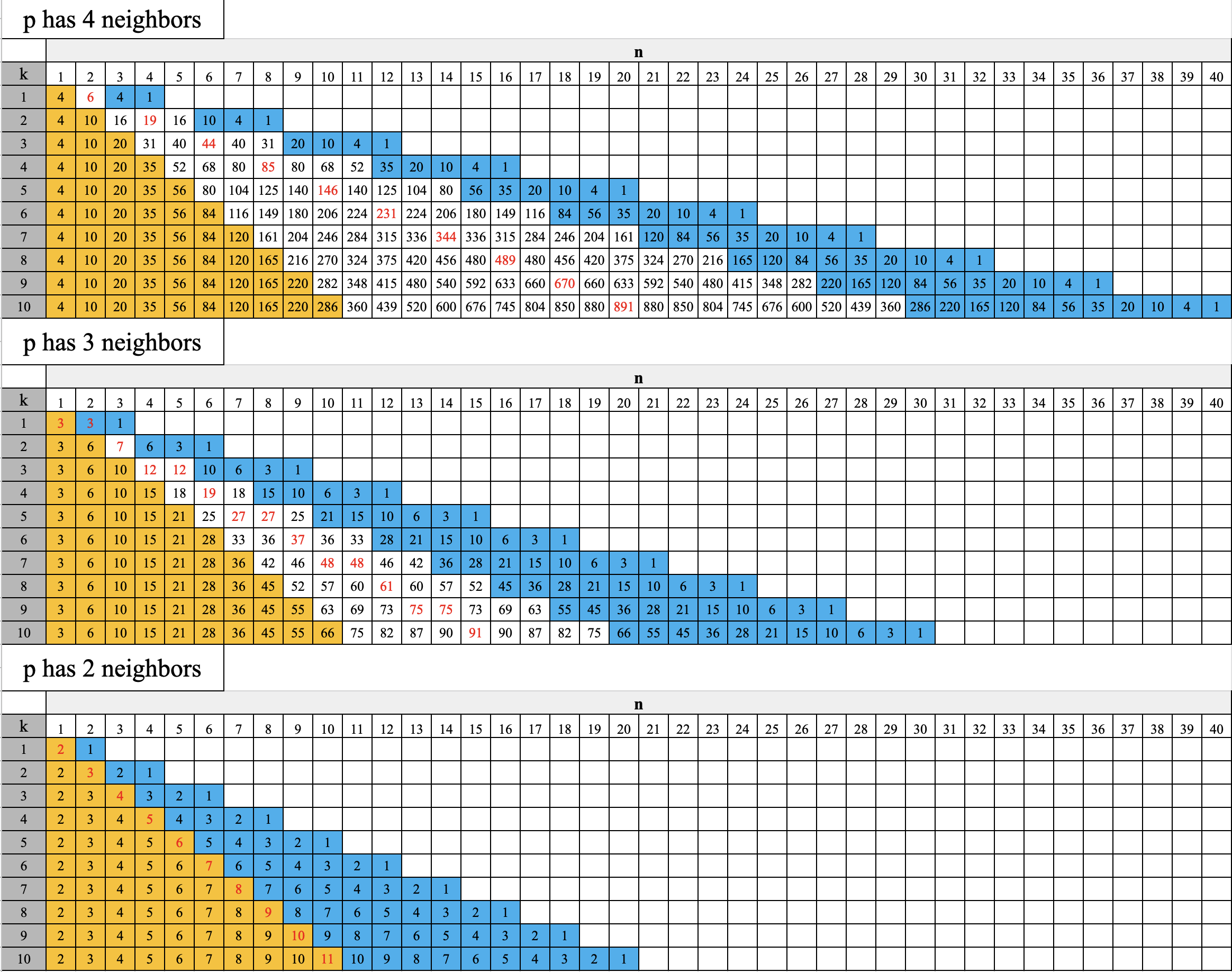}
    \caption{Entries in each of the three tables indicate the number of connection configurations that are possible for a node of magnitude $n$ with 4, 3, and 2 neighbors and the indicated $k$. Values in orange indicate when the number of configurations ceases to increase with $k$. Values in blue show what values appear again in the tail of the row below. Values in red are the maximums in each row.  }
    \label{fig:where2start_table}
\end{figure}

In Figure \ref{fig:where2start_table} we give a table that shows how many connection configurations exist given a node of magnitude $n$, the value for $k$, and the number of neighbors the node has. 
We note some observations based on this table. 
\begin{itemize}
    \item For a fixed $k$ (a row in the table), the number of connection configurations is symmetric as the magnitude $n$ of a node increases. The reader is encouraged to figure out why! Here is a hint: what similarities exist between the connection configurations for the node of magnitude 7 with 4 neighbors and $k=2$ (example above) and a node of magnitude 1 with 4 neighbors and $k=2$?  
    \item The maximum number of connection configurations possible for a fixed $k$ is achieved at magnitude $n = \lfloor\frac{rk}{2}\rfloor$ where $r$ is the number of neighbors  the node has. 
    \item The "tail" of each row in the table for a node with 4 neighbors can be generated by the following formula \cite{sloane2003line}:
    $$ \frac{i(i+1)(i+2)}{6}\hspace{.2in} \text{ for } \hspace{.2in} i =1, 2, 3, \dots$$
    \item For 3 neighbors, the "tail" of each row is simply the triangular numbers:  
    $$ \frac{i(i+1)}{2}\hspace{.2in} \text{ for } \hspace{.2in} i =1, 2, 3, \dots.$$
    \item \textbf{Where do we start?} The table gives us an idea of where to start. The smaller numbers in the table gives us a higher chance of success in creating connections based on $\omega^*$. The table also shows us what nodes are the low-hanging fruit. For example, we see that if a node has $r$ neighbors and has magnitude $magn(p) = rk$, then there is only 1 way to create the needed connections. 
\end{itemize}

\subsection{The $\tau$ algorithm}
We now describe an algorithm which creates the guaranteed connections in a $k$-grid. We also show that if this algorithm solves a $k$-grid, then the solution is unique. 

To begin we will formalize the idea of creating connections between nodes and how this creates a new $k$-grid. 

\begin{defn}
For any $k$-grid $\Gamma_k$, define the \textbf{builder function} for node $p \in \Gamma_k$ by,
$$u\left(p,\omega^{*}(p)\right) = \Gamma_k^{'}$$
Where $\Gamma_k^{'}$ is $\Gamma_k$ with the following changes:
\begin{enumerate}
    \item $p \in \Gamma_k^{'}$ where $p = \big[x_p, y_p, n_p - ||\omega^{*}||\big]$
    \item If $q \in \Gamma_k^{'}$ is the $l_p$ neighbor of $p$, then $q = \big[x_q, y_q, n_q - count(l,\omega^{*}) \big]$.
    \item For each $l \in W$ and pair of nodes $p$ and $l_p$, if $count(l,\omega^{*})>0$, we construct $count(l,\omega^{*})$ connections between $p$ and $l_p$.
\end{enumerate}
\end{defn}

At their core, builder functions build connections between a node and its neighbors, reduce the magnitude of the node and reduce the magnitude of the neighboring nodes. We will use these functions to build up a solution (should one exist) for a general $k$-grid. 
We would like to apply the builder functions to each node in a $k$-grid again and again in hopes of finding additional guaranteed connections. 
We now describe the algorithm which computes and draws the guaranteed connections between nodes. This algorithm will continue until the grid is solved or all guaranteed connections have been exhausted.

$$\textbf{Tau Algorithm}$$
\begin{enumerate}
    \item Choose a node $p$ with $r$ neighbors with $magn(p)=rk$. In this case $k$ connections must be formed with each neighbor. Draw the needed connections for each node $p$ with this property. 
    \item Next, choose a node with one neighbor. Again there is only one destination for those connections! Complete this step for every node $p$ with one neighboring node.
    \item Now choose a node with one remaining incomplete neighbor. With all other neighbors completed, there is only one place left to draw connecting lines. Draw the needed connections for each node $p$ with one incomplete neighbor.
    \item 
    If there are no remaining nodes with one incomplete neighbor, choose a node $p$ via the following guidelines which are based on the intuition we gather from Figure \ref{fig:where2start_table}. Then compute and draw the connection to create $u(p,\omega^*)$.
        \begin{enumerate}
            \item Start with the nodes which have the smallest number of neighbors. Based on Figure \ref{fig:where2start_table}, we know these typically have the fewest possible connection configurations.
            \item Of these nodes, start with those with the largest and smallest magnitudes (putting off magnitudes that are near $\lfloor \frac{rk}{2} \rfloor$ where $r$ is the number of neighbors that $p$ has).
        \end{enumerate}
    \item Repeat Step 3-4 until either no nodes with magnitude larger than 0 remain or $\omega^*(p) = \emptyset$ for each remaining node with non-zero magnitude. 
\end{enumerate}
For any $k$-grid $\Gamma_k$, we denote the resulting grid from the $\tau$ algorithm as $\tau(\Gamma_k)$.

\noindent \textit{\textbf{Theorem }}
If $\tau(\Gamma_k)$ solves $\Gamma_k$, then $\Gamma_k$ is uniquely solved. \\

\noindent \textit{\textbf{Proof. }}
Note that by construction, the word $\omega^*$ for a node $p$ gives a \emph{subset} of the connections that will exist in any solution to the grid which contains $p$. At each step in the $\tau$ algorithm, we calculate the $\omega^*$ for a node and create the connections indicated by the $\omega^*$. This means that if $\tau(\Gamma_k)$ is a solved $k$-grid, then the set of connections indicated by repeated computations of $\omega^*$ gives the only configuration of connections.  $\blacksquare$
\medskip

\begin{figure}
    \centering
    \includegraphics[scale = .4]{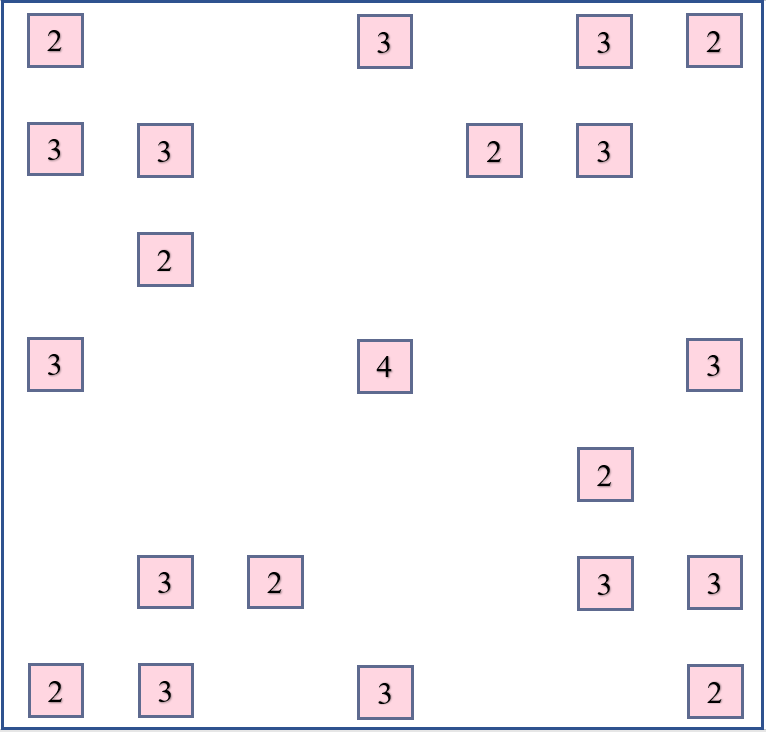}
    \caption{This $k$-grid which we call the Pinwheel has a unique solution and yet none of it's connections can be found using the $\tau$ algorithm.  }
    \label{fig:Pinwheel}
\end{figure}

This theorem tells us that if the $\tau$ algorithm solves a $k$-grid, then we can be sure that we have found the only solution possible. One may wonder if there is hope that the above theorem be bidirectional. Unfortunately, it is not true that if a $k$-grid has a unique solution, then the $\tau$ will be able to solve it.  Consider the $k$-grid from Figure \ref{fig:Pinwheel}. When $k=2$, this $k$-grid which we call the Pinwheel is uniquely solved. 
Every node has $r>1$ neighbors with $magn(p)<rk$, meaning no connections will be drawn from the first three steps. When computing $\omega^*$ for any node we get $\omega^*(p) = \emptyset$. So $\tau$ is unable to make any connections in the Pinwheel! 

\begin{figure}
    \centering
    \includegraphics[scale = .45]{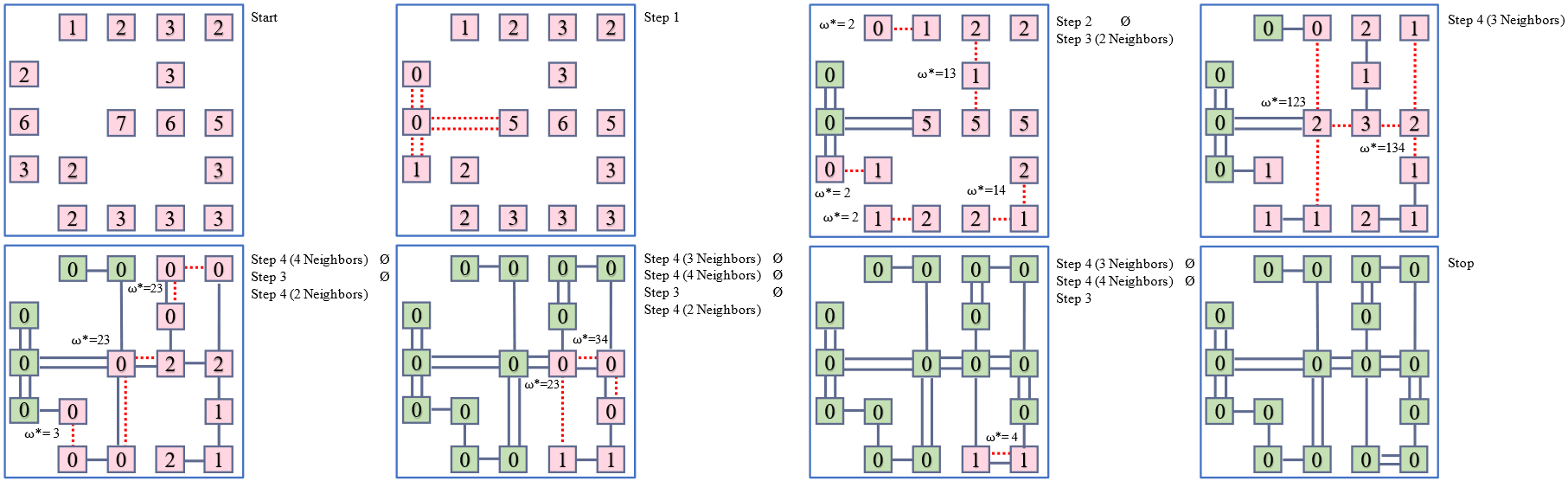}
    \caption{A solved $k$-grid whose connections where determined by $\tau$. This $k$-grid is also included in Numberlink and listed as Expert \cite{pzzldom}.}
    \label{fig:TauEx}
\end{figure}

\section{Future Work}

In this work we have created a framework for working with Numberlink puzzles as graphs with certain conditions on the connections between graphs. We have also described an algorithm for computing the guaranteed connections of a node and have shown that in the case that our algorithm solves a $k$-grid, we have found the one unique solution. There are many possible generalizations to this framework and open questions. Below we list a few of these. 
\begin{itemize}
    \item The choice of a square to represent each node is arbitrary and used to emphasize the restriction on the number of neighbors of a node (maximum 4). One could look at graphs whose nodes are polygons (such as an octagon) where connections can be drawn to $m$ neighbors. 
    \item As mentioned in section 3, some nodes have fewer options for how connections can be configured based on $k$ and the number of neighbors the node has. It would be helpful to find a formula to express the number of ways that $n$ connections can be configured for a node based on $k$ and the number of neighbors the node has. This can then tell us what nodes are good starting points for solving the puzzles. \item Having a formula as described above may lead to a more refined algorithm for $\tau$ and may give a way of calculating how many steps will be needed in order to solve a $k$-grid. This can then help find the solution with the least number of steps. 
    \item We know that many $k$-grids are solvable regardless of the value of $k$ and that if a grid is solvable for $k$ then it is also solvable for $k' > k$. Can we find a process to easily compute the smallest value of $k$ such that the grid is solvable? 
    \item In what way can this framework be used for applications in areas like operations research and optimization. 
\end{itemize}


\bibliographystyle{unsrt}
\bibliography{references}

\end{document}